%% file: correctedversion.tex
\documentclass[11pt]{llncs}
\usepackage{amsmath, amssymb}
\usepackage{fontenc,graphicx}
\usepackage{url}
\usepackage{enumerate}
\usepackage[all]{xy}

\newcommand{\sS}{\mathcal{S}}
\newcommand{\PP}{\mathbb{P}}
\newcommand{\CC}{\mathbb{C}}
\newcommand{\FF}{\mathbb{F}}
\newcommand{\AAb}{\mathbb{A}}
\newcommand{\ZZ}{\mathbb{Z}}
\newcommand{\Mm}{\mathsf{M}}

\DeclareMathOperator{\Aut}{Aut}

\DeclareMathOperator{\echr}{char}

\DeclareMathOperator{\End}{End}
\DeclareMathOperator{\Gal}{Gal}

\DeclareMathOperator{\Hom}{Hom}

\DeclareMathOperator{\Jac}{Jac}

\DeclareMathOperator{\Sym}{Sym}

\institute{ \'Equipe GAATI,
         Universit\'e de la Polyn\'esie Fran\c{c}aise, 
         BP 6570 - 98702 Faa'a - Tahiti - Polyn\'esie fran\c{c}aise \\
\email{roger.oyono@upf.pf }
\and Institut de Math{\'e}matiques de Luminy,
         UMR 6206 du CNRS,
         Luminy, Case 907, 13288 Marseille, France \\
\email{ritzenth@iml.univ-mrs.fr}
}

\author{Roger Oyono\inst{1}  \and Christophe Ritzenthaler\inst{2} 
\thanks{The second author
acknowledges the financial support of the grant MTM2006-11391 from the Spanish MEC and of the grant ANR-09-BLAN-0020-01 from the French ANR}
}
\title{On rationality of the intersection points of a line with a plane quartic}

\begin{document}
\maketitle

\begin{abstract}
We study the rationality of the intersection points of certain lines and  smooth plane quartics $C$ defined over $\FF_q$. For $q\geq 127$, we prove the existence of a line such that the intersection points with  $C$ are all rational. Using another approach, we further prove the existence of a tangent line with the same property as soon as   $q \geq 66^2+1$. Finally, we study the probability of the existence of a rational flex on $C$ and exhibit a curious behavior when $\echr \FF_q=3$.

\medskip

\keywordname{ smooth plane quartics, rationality, intersection points, tangent line, flex.}

\subclassname{ 11G20, 14G05, 14G15, 14H45, 14N10.}

\end{abstract}

\input intro2

\input version8

\input Flexes_char_3

\bibliographystyle{plain}
\bibliography{synthbib}

\end{document}

%% file: intro2.tex
\section{Introduction}

In computational arithmetic geometry, it is an important task to develop an efficient group law for the Jacobian variety of algebraic curves defined over a finite field. One of the most important and recent application of such efficient arithmetic comes from cryptography \cite{Koblitz,Koblitz2,Cantor}. In \cite{FOR}, the authors introduced an efficient algorithm to perform arithmetic in the Jacobian of smooth plane quartics. The  presented algorithm depends on the existence of a rational line $l$ intersecting the quartic in rational points only. Moreover, the more special $l$ is (for instance, tangent, tangent at a flex,\ldots), the better is the complexity of the algorithm. Motivated by the above efficiency argument, we prove here the following theorems.

\begin{theorem}\label{conditionetoile1}
Let $C$ be a smooth plane  quartic over the finite field $\FF_q$ with $q$ elements. If $q \geq 127,$ then there
exists a line $l$ which intersects $C$ at rational points only.
\end{theorem}


\begin{theorem}\label{conditionetoile2}
Let $C$ be a smooth plane quartic over $\FF_q$. If $q \geq 66^2+1,$ then there
exists a tangent to $C$ which intersects $C$ at rational points only.
\end{theorem}

In  \cite{FOR}, the authors gave heuristic arguments and computational evidences that the probability for a plane smooth quartic over a finite field to have a rational flex is about $0.63$. In this article, we present a `proof' which depends on a conjectural analogue over finite fields of a result of \cite{harris-enu} on the Galois group $G$ of the $24$ flexes of a general quartic. Harris proved that, over $\CC$, the group $G$ is as big as possible, namely the symmetric group $S_{24}$. 
Unfortunately, Harris's proof uses monodromy arguments which cannot be adapted so easily in positive characteristic. Worse, we surprisingly found that Harris's result is not valid over fields of characteristic $3$ and we prove that for this field the Galois group is $S_8$. This is a consequence of the peculiar fact that, in characteristic $3$, a smooth plane quartic $C$ has generically only $8$ flexes (with multiplicity $3$), which belong to a conic. We suspect that this is the only exceptional case.\\

The methods used for these three problems are various and can be  generalized or adapted to other questions. This was our principal motivation to write down our approaches in the case of quartics. It illustrates also the very unusual behavior of special points and lines in small characteristics. \\

Coming back to our initial motivation, it appears nowadays unlikely (due to recent progress in index calculus \cite{DiTh}) that smooth plane quartics may be used for building discrete logarithm cryptosystems. 
 However, it is interesting to mention that the complexity analysis of the index calculus attack of \cite{DiTh} uses an asymptotic bound for the number of  lines intersecting a smooth plane quartic in four distinct rational points, in the spirit of Theorem \ref{conditionetoile1}.\\

 The paper is organized as follows: Section \ref{wfact} gives a brief overview on the possible geometric intersections of a line and a smooth plane quartic. 
  In Section \ref{tcase}, we give a proof of Theorem \ref{conditionetoile1} using Chebotarev density theorem for covers of curves. In Section \ref{ttcase}, we prove Theorem \ref{conditionetoile2} using the tangential correspondence and its associated curve $X_C \subset C \times C$. When $\echr \FF_q \ne 2$, the crucial point is to prove that $X_C$ is geometrically irreducible in order to apply Hasse-Weil bound for (possibly singular) geometrically irreducible curves.
 For the characteristic $2$ case, $X_C$ is not always absolutely irreducible. However, if it is not, we prove it is reducible (over $\FF_q$) and we can apply Hasse-Weil bound on each component.
Finally, in Section \ref{open}, we also address the question of the probability of existence of a rational flex and show that in characteristic $3$, the flexes are on a conic.\\

\noindent {\bf Conventions and notation.} In the following, we denote by
$(x:y:z)$  the  coordinates in $\mathbb{P}^2$, and by $(x,y)$ the coordinates in $\mathbb{A}^2$. Let $p$ be a prime or $0$ and $n \geq 1$ an integer. We use the letter $K$ for an arbitrary field of characteristic $p$ and let $k=\FF_q$ be a finite field with $p \ne 0$ and $q=p^n$ elements. When $C$ is a smooth geometrically irreducible projective curve, we denote by $\kappa_C$ its canonical divisor. Operators such as $\Hom$, $\End$ or $\Aut$ applied to varieties
over a field $K$ will always refer to $K$-rational homomorphisms and endomorphisms.\\

\noindent {\bf Acknowledgments.} We would like to thank Noam Elkies and Pierre D\`ebes for their suggestions and references in Section \ref{flexcar3}. We are also very grateful to Felipe Voloch for showing us how to prove Theorem \ref{conditionetoile2} in characteristic $2$.


%% file: version8.tex
\section{Structure of the canonical divisor}\label{wfact}

In this section, we recall some geometric facts about special points and lines on a plane smooth quartic. Let $K$ be an algebraically closed field of characteristic $p$ and   $C$ be a smooth (projective) plane quartic defined over $K$. The curve $C$ is a non hyperelliptic genus $3$ curve which is canonically embedded.
Hence the intersection of $C$ with a line $l$ are the positive canonical divisors of $C$. There are $5$ possibilities for the intersection
divisor  of $l$ and $C$ denoted $(l\cdot C)=P_1 + P_2 + P_3 + P_4$:
\begin{enumerate}[\textrm{case} 1]
\item \label{item:wfact1} The four points
  are pairwise distinct. This is the generic position. 
\item \label{item:wfact2}   $P_1=P_2$, then $l$ is tangent to $C$ at $P_1$.  
\item \label{item:wfact3}   $P_1=P_2=P_3$. The point $P_1$ is then called a
  \emph{flex}. As a linear intersection also represents the canonical
  divisor $\kappa_C$, these points are exactly the ones where a regular
  differential has a zero of order $3$. 
  The curve $C$ has infinitely many flexes if and only if $ p=3$ and
$C$ is isomorphic to $x^4+y^3z+yz^3=0$ which is also isomorphic to 
the Fermat quartic $x^4+y^4+z^4=0$ and to the Klein quartic $x^3y+y^3z+z^3x=0$. This is a \emph{funny curve} in the sense of \cite[Ex.IV.2.4]{Hart} or a \emph{non classical} curve in the sense of  \cite[p.28]{Torres}.
 On the contrary, if $C$ has finitely many flexes, then these points are
the Weierstrass points of $C$ and the sum of their weights is $24$.
\item \label{item:wfact4}  $P_1=P_2$ and $P_3=P_4$. The line $l$ is called a
 \emph{bitangent} of the curve $C$ and the points $P_i$ \emph{bitangency points}.
  If
  $p \ne 2$,  then $C$ has exactly $28$ bitangents (see for instance \cite[Sec.3.3.1]{phdritz}). If
  $p  = 2$, then $C$ has respectively $7,4,2$, or $1$ bitangents,
  if the $2$-rank of its Jacobian is respectively $3,2,1$ or $0$ \cite{stvo}. Recall that the $p$-rank 
 $ \gamma$ of an abelian variety $A/K$ is defined by $\# A[p](K)=p^{\gamma}$.
\item \label{item:wfact5} $P_1=P_2=P_3=P_4$. The point $P_1$ is called a
  \emph{hyperflex}. Generically, such a point does not exist. More precisely, the locus of quartics
  with at least one hyperflex is of codimension one in the moduli space $\Mm_3$ (see \cite[Prop.4.9,p.29]{vermeulen}). 
If $p=3$ and $C$ is  isomorphic to the Fermat quartic then the number of hyperflexes of $C$ is equal to $28$. Otherwise, the weight of a hyperflex is greater than or equal to $2$, so there are less than $12$ of them \cite{stvo2}. Note that the weight of a hyperflex is exactly $2$ when $p>3$ or $0$ and that the weight of a flex which is not a hyperflex is $1$ if $p \ne 3$. See also \cite{viana} for precisions when $p=2$ and Section \ref{flexcar3} when $p=3$. 
\end{enumerate}

A point $P$ can even be more special. Let $P \in C$ be a point and let us denote $\phi_P :  C \to |\kappa_C-P|=\PP^1$ the degree three map induced by the linear system $|\kappa_C-P|$. If this cover is Galois, such a point $P$ is called a (inner)
\emph{Galois point} and we denote $\Gal(C)$ the set of Galois points of $C$. 
\begin{lemma}[\cite{fukasawa1,fukasawa2,fukasawa3,fukasawa4}] \label{galoispoints}
Let $C$ be a smooth plane quartic defined over $K$. The number of Galois points is at most $4$ if
$p  \ne 3$ and at most $28$ if $p =3$. Moreover, the above bounds are reached, respectively by the curve $yz^3+x^4+z^4=0$ and the Fermat quartic.
\end{lemma}

\noindent In the sequel, we will need the first item of the following lemma.

\begin{lemma} \label{classification}
Let $C$ be a smooth plane quartic defined over  $K$. There is always a bitangency point which is not a
hyperflex unless
\begin{itemize}
\item  $p  = 3$ and $C$ is   isomorphic to the Fermat
 quartic,
\item  $p =2$ and $\Jac C$ is supersingular,
\item $p =2$ and $C$ is isomorphic to a $2$-rank one quartic
 $$(a x^2+ b  y^2+ c z^2+ d
  xy)^2+xy(y^2+xz)=0$$ with $ac \ne 0$,
\item $p  =2$ and $C$ is isomorphic
to a $2$-rank two quartic $$(a x^2+ b  y^2+ c z^2)^2+xyz(y+z)=0$$
with $abc \ne 0$ and $b+c \ne 0$.
\end{itemize}
\end{lemma}
\begin{proof}
According to Section \ref{wfact} case \ref{item:wfact5}, the number of
hyperflexes when $C$ is not isomorphic to the Fermat quartic and
$p =3$, is less than $12$. On the other hand, if $p \ne 2$, a
curve $C$ has $ 28$ bitangents, and thus there are at least $28 \cdot 2-12 \cdot 2=32$ bitangency points which are
not hyperflexes. However, for the Fermat quartic
in characteristic $3$, all bitangency points are hyperflexes.\\
There remains to look at the case $p=2$ for which we use the
classification of \cite{wa}, \cite{NR06}. The quartic $C$ falls
into four categories according to its number of bitangents:
\begin{enumerate}
\item if $C$ has only one bitangent, then $C$ is isomorphic
  to a model of the form $Q^2=x(y^3+x^2z)$ where $Q=a x^2+b y^2+c z^2+ d x y + e y
  z+f z x$ and $c \ne 0$. The unique bitangent $x=0$ intersects  $C$ at points $(x:y:z)$ satisfying $b
  y^2+ c z^2 + e y z=0\, .$ Therefore, $C$ has a  hyperflex if and only if $e=0$,
  \textit{i.e.} $C$ falls into the subfamily $\sS
 $ \cite[p.468]{NR06} of
  curves whose Jacobian is supersingular. Conversely, any curve in $\sS$ has a hyperflex.
\item if $C$ has two bitangents, then $C$ is isomorphic  to a
  model of the form $Q^2=xy (y^2+x z)$ where $Q=a x^2+b y^2+c z^2+ d x y + e y
  z+f z x$ and $a c \ne 0$. All bitangency points are then hyperflexes
  if and only if $e=f=0$.
\item if $C$ has four bitangents, then $C$ is isomorphic
  to $Q^2=xyz(y+z)$ with  $$Q=a x^2+b y^2+c z^2+ d x y + e y
  z+f z x \; \textrm{and} \; abc \ne 0, b+c+e \ne 0.$$ All bitangency points are
  hyperflexes if and only if $d=e=f=0$.
\item if $C$ has seven bitangents, then $C$ is isomorphic
  to  $Q^2=xyz(x+y+z)$ with  $Q=a x^2+b y^2+c z^2+ d x y + e y
  z+f z x$ and some open conditions on the coefficients
  \cite[p.445]{NR06}. The bitangents are $$\{x,y,z,x+y+z,x+y,y+z,x+z\}.$$
 Suppose that the intersection points of $x=0$, $y=0$ and
  $z=0$ with $C$ are hyperflexes, then $d=e=f=0$. Moreover $x+y=0$
  gives a hyperflex if and only if $(a+b) y^2+ yz+c z^2$ is a perfect
  square, which is never possible. \qed
\end{enumerate} 
\end{proof}

%


\section[Proof]{Proof of Theorem \ref{conditionetoile1}}\label{tcase}

\noindent Let $q \geq 127$ be a prime power. Note that, as an easy consequence of Serre-Weil bound, we know that $$\#C(k) \geq q+1-3 \cdot \lfloor 2 \sqrt{q} \rfloor =62.$$ For the proof, we follow the same strategy as \cite[p.604]{DiTh}. The lines  intersecting $C$ at $P$
are in bijection with the divisors in the complete linear system
$|\kappa_C-P|$. We wish to estimate the number of completely split
divisors in this linear system, since such a divisor defines a line solution of Theorem \ref{conditionetoile1}. To get the existence of such a divisor, we will use an effective Chebotarev density
theorem for function fields, as in \cite[Th.1]{murty}.
This theorem assumes that the cover is Galois  but we can reduce to this case thanks to the following lemma.

\begin{lemma} \label{galoisex}
Let $K/F$ be a finite separable  extension of function fields over a
finite field. Let $L$ be the Galois closure of $K/F$. A
place of $F$ splits completely in $K$ if and only if it splits
completely in $L$.
\end{lemma}
\begin{proof}
It is clear that, if a place $P \in F$ splits completely in $L$,
it splits completely  in $K$. Conversely, let $G$ be the Galois
group of $L/F$ and $H$ be the Galois subgroup of $L/K$. By
construction (see \cite[A.V.p.54]{bourbaki}), $L$ is the
compositum of the conjugates $K^{\sigma}$ with $\sigma \in G/H$. If
a place $P \in F$ splits completely in $K$, it splits completely in
each of the $K^{\sigma}$. It is then enough to apply
\cite[Cor.III.8.4]{stich} to conclude.   \qed
\end{proof}

 We
consider the separable geometric cover $\phi_P : C \to
| \kappa_C-P|=\PP^1$ of degree $3$ induced by the linear system $| \kappa_C-P|$.
We may assume that no rational point in $\PP^1$ is ramified for $\phi_P$. Otherwise, it is easy to see that the fiber of $\phi_P$ above this point has only rational points and the line defined by these points intersects the quartic in rational points only. Theorem \cite[Th.1]{murty} boils down to the following proposition.


\begin{proposition} \label{chebobo}
\begin{enumerate}[a)]
\item If the cover $\phi_P$ has a non-trivial automorphism, 
then the number $N$ of completely split divisors in $| \kappa_C-P|$
satisfies
$$\left|N-\frac{q+1}{3}\right| \leq 2 \sqrt{q}+ |D|$$
where $|D|=\sum_{y \in \PP^1, \textrm{ramified}} \deg y$. \label{caseok}
\item If the cover $\phi_P$ has a non-trivial $\bar{k}$-automorphism not defined
  over $k$, then there are no completely split divisors in $| \kappa_C-P|$. \label{badcase}
\item  If the cover $\phi_P$ has no non-trivial $\bar{k}$-automorphism, then the
  number $N$ of completely split divisors in $| \kappa_C-P|$ satisfies
$$\left|N-\frac{q+1}{6}\right| \leq  \sqrt{q}+ |D|.$$
\end{enumerate}
\end{proposition}





\begin{proof}[of Theorem \ref{conditionetoile1}]
\noindent To avoid case \eqref{badcase} of Proposition \ref{chebobo}, it is enough that $P$ is not a Galois point (and we will avoid case \eqref{caseok} as well). By Lemma \ref{galoispoints},  we know that the number of such $P$  is less than $28$. So let the point $P \in C(k) \setminus (\Gal(C) \cap C(k))$. Then the cover $\phi_P$ has a completely split divisor if
\begin{equation} \label{cheboo}
\frac{q+1}{6}  >  \sqrt{q}+ |D|.
\end{equation}
Using Riemann-Hurwitz formula we get
$$|D| \leq (2 \cdot 3-2) - 3 \cdot (0-2)=10.$$ 
Hence the inequality \eqref{cheboo} is satisfied as soon as $q \geq 127$. \qed
\end{proof}

\begin{remark}
We do not pretend that our lower bound $127$ is optimal. In \cite{Ballico}, the converse problem is considered (\textit{i.e.} the existence of a plane  (not necessarily smooth) curve with no line solution of Theorem \ref{conditionetoile1}) but their bound, $3$, is also far from being optimal in the case of quartics. Indeed, by \cite{HLT}, for $q=32$, there still exists a pointless smooth plane quartic for which of course there is no line satisfying Theorem \ref{conditionetoile1}.
\end{remark}
%




\section[Proof]{Proof of Theorem \ref{conditionetoile2}} \label{ttcase}

Let $C$ be a smooth plane quartic defined over a field $K$. Let the map $T : C \to \Sym^2(C)$ be the {\em tangential correspondence} which sends a point $P$ of $C$ to the divisor $
(T_P(C) \cdot C) - 2 P$.  We associate
to $T$ its correspondence curve $$X_C=\{(P,Q) \in C \times C : Q \in
T(P)\}$$ which is defined over $K$. Our goal is to show that when $K=k=\FF_q$ with $q >66^2$, then $X_C$ has a rational point, \textit{i.e.} there is
$(P,Q) \in C(\FF_q)^2$ such that $(T_P(C) \cdot C) = 2 P +Q+R$ for some point $R$, necessarily in $C(\FF_q)$.
 Thus, the
tangent $T_P(C)$ is a solution of Theorem \ref{conditionetoile2}.

\subsection{Study of the geometric cover $X_C \to C$}

 Let $\pi_i : X_C \to C$, $i=1,2,$ be the
projections on the first and second factor. The morphism $\pi_1$ is a $2$-cover between these two projective curves.

\begin{lemma}\label{regpoints}
Let $K$ be an algebraically closed field of characteristic $p$. The projection  $\pi_1 : X_C \longrightarrow C$ has the following properties:
\begin{enumerate}
\item The ramification points of $\pi_1$
 are the bitangency points of $C$,
\item $\pi_1$ is separable,
\item The point $(P,Q) \in X_C$ such that $P,Q$ are bitangency points
and $P$ is not a hyperflex (\textit{i.e.} $P \ne Q$) is a regular point if
and only if $p \ne 2$,
\item If $p \ne 2$, the only possible singular points of $X_C$
are the points $(P,P)$ where $P$ is a hyperflex of $C$.
\end{enumerate}
\end{lemma}
\begin{proof}
\noindent The first property is an immediate consequence of the definition of a
bitangent.

\noindent If $\pi_1$ is not separable then $p =2$ and $\pi_1$ is purely
inseparable. Thus $\# \pi_1^{-1}(P)=1$ for all $P \in C$, \textit{i.e.} all
$P$ are bitangency points. This is impossible since the number of bitangents is finite
(less than or equal to $7$).

\noindent Let $F(x,y,z)=0$ be an equation of  $C$. Let $Q \ne P$ be a point of $C$ defining a point $(P,Q)$  in $X_C \setminus \Delta$
where $\Delta$ is the diagonal of $C \times C$. For such points, it
is easy to write local equations as follows. We can suppose that
 $P=(0:0:1)=(0,0)$, $Q=(1:0:1)=(1,0)$ and assume that
  $f(x,y)=F(x,y,1)=0$ is an equation of the affine part of $C$. Then,
  if we consider the curve $Y_C$ in $\AAb^4(x,y,z,t)$ defined by
$$\begin{cases}
f(x,y) &= 0\,, \\
f(z,t) &= 0\, , \\
\frac{\partial f}{\partial x}(x,y) (z-x)+ \frac{\partial f}{\partial
y}(x,y)
(t-y) &= 0\, , 
\end{cases}$$  $Y_C \setminus \Delta$ is
an open subvariety of  $X_C$ containing $(P,Q)$. The Jacobian matrix at the
point $(P,Q)=((0,0),(1,0))$ is equal to  {\tiny
$$\left(\begin{array}{ccc} \frac{\partial f}{\partial x}(0,0) & 0
& 
  \frac{\partial^2
  f}{\partial x^2}(0,0)-\frac{\partial f}{\partial x}(0,0) \\
\frac{\partial f}{\partial y}(0,0)  & 0 & 
   \frac{\partial^2
  f}{\partial x \partial y }(0,0) -\frac{\partial f}{\partial y}(0,0)  \\
0 & \frac{\partial f}{\partial x}(1,0) & \frac{\partial f}{\partial x}(0,0) \\
0 & \frac{\partial f}{\partial y}(1,0)  & \frac{\partial f}{\partial
y}(0,0)
\end{array}\right).$$}
\noindent Now, if $P$ and $Q$ are bitangency points, then the tangent at these points
is $y=0$, so $\frac{\partial f}{\partial x}(0,0)=\frac{\partial
f}{\partial x}(1,0)=0$. The only non-trivially zero minor
determinant of the matrix is then
 $$\frac{\partial f}{\partial y}(1,0) \cdot \frac{\partial
f}{\partial y}(0,0) \cdot \frac{\partial^2 f}{\partial x^2}(0,0).$$
So $(P,Q) \in X_C$ is not singular if and only if $\frac{\partial^2
f}{\partial x^2}(0,0) \ne 0$. This can never be the case if
$p =2$, so we now suppose that $p  \ne 2$. We can always
assume that the point $(0:1:0) \notin C$ and we write
$$f(x,y)=x^4+x^3 h_1(y)+x^2 h_2(y)+x h_3(y)+h_4(y),$$
where $h_i$ are polynomials (in one variable) over $K$ of degree $\leq i$. Since $y=0$
is a bitangent at $P$ and $Q$, we have
$$f(x,0)=x^2(x-1)^2=x^4-2x^3+x^2,$$
and thus $h_2(0)=1$. Now $$\frac{\partial^2 f}{\partial x^2}(0,0)=2
h_2(0) \ne 0.$$ Finally if $(P,Q) \in X_C$ is not ramified for $\pi_1$, it is a
smooth point. This proves the last assertion. \qed
\end{proof}

We want to apply the following version of Hasse-Weil bound to the curve $X_C$.
\begin{proposition}[\cite{aubry}] \label{paubry}
Let $X$  be a geometrically irreducible curve of
arithmetic genus $\pi_X$ defined over $\FF_q$. Then
$$|\#X(\FF_q) - (q+1)| \leq 2\pi_X \sqrt{q}.$$
In particular if $q \geq (2 \pi_X)^2$ then $X$ has a rational point.
\end{proposition}

Hence, to finish the proof, we need to show that $X_C$ is geometrically irreducible and then  to compute its arithmetic genus. For the first point, we use the
following easy lemma for which we could not find a reference.

\begin{lemma} \label{irredu}
Let $\phi : X \to Y$ be a separable morphism of degree $2$ between two
projective curves defined over an algebraically closed field $K$
such that
\begin{enumerate}
\item $Y$ is smooth and irreducible,
\item there exists a point $P_0 \in Y$ such that $\phi$ is ramified at
  $P_0$ and $\phi^{-1}(P_0)$ is not singular.
\end{enumerate}
Then $X$ is irreducible.
\end{lemma}
\begin{proof}
Let $s: \tilde{X} \to X$ be the normalization of $X$ and
$\tilde{\phi}=\phi \circ s : \tilde{X} \to Y$. Due to the second
hypothesis, $\tilde{\phi} : \tilde{X} \to Y$ is a separable,
ramified $2$-cover. Clearly, $X$ is geometrically irreducible if and
only if
$\tilde{X}$ is.\\
Let us assume that $\tilde{X}$ is not  irreducible. There exist smooth projective curves
 $\tilde{X}_1$ and $\tilde{X}_2$  such
that $\tilde{X}= \tilde{X}_1 \cup \tilde{X}_2$. Then, consider for
$i=1,2$, $\tilde{\phi}_i= \tilde{\phi}_{|\tilde{X}_i} : \tilde{X}_i
\to Y$. Each of these morphisms is of degree $1$ and since the
curves are projective and smooth, they define an isomorphism between
$\tilde{X}_i$ and $Y$.\\
Since $P_0 \in Y$ is a ramified point,
$\tilde{\phi}_1^{-1}(P_0)= \tilde{\phi}_2^{-1}(P_0)$.  It follows that
$$\tilde{\phi}^{-1}(P_0) \in \tilde{X}_1 \cap \tilde{X}_2$$ so that $\tilde{\phi}^{-1}(P_0)$ is singular, which contradicts the hypothesis. \qed
\end{proof}

\subsection{Proof of Theorem \ref{conditionetoile2}: the case $p \ne 2$}
\label{sub:proof}
Let us first start with $C$ $\bar{k}$-isomorphic to the Fermat quartic in characteristic $3$. By Section  \ref{wfact} case \ref{item:wfact3}, all its points are
flexes. So if there exists $P \in C(k)$,  then the tangent at $P$ cuts
$C$ at $P$ and at another unique point which is again rational over $k$.  Now when $q>23$ and $q\neq 29, 32$, it is proved in
\cite{HLT} that a genus $3$ non-hyperelliptic curve over $\FF_q$ has always a
rational point and the result follows.\\
We suppose  that $C$ is a  smooth plane quartic not $\bar{k}$-isomorphic to the Fermat quartic if $p=3$. As we assumed that  $p \neq 2$,  by Lemma \ref{irredu}, Lemma \ref{classification} and Lemma 
\ref{regpoints}, we conclude that $X_C$ is  an  geometrically irreducible projective curve. Moreover, if we
assume that $C$ has no hyperflex, then $X_C$ is smooth and we can compute
its genus $g_{X_C}$ using Riemann-Hurwitz formula for the $2$-cover $ \pi_1 : X_C \to C$ ramified
over the $2 \cdot 28$ bitangency points. In fact,
$$ 2g_{X_C} -2= 2(2\cdot 3-2) +56\, ,
$$
and thus $g_{X_C}=33$. The family of curves $X_C$ is flat over the locus of smooth plane quartics $C$ by \cite[Prop.II.32]{eisenbud-harris}. As the
arithmetic genus $\pi_{X_C}$ is constant in flat families \cite[Cor.III.9.10]{Hart} and equal to $g_{X_C}$ for smooth $X_C$ \cite[Prop.IV.1.1]{Hart}, we get that
$\pi_{X_C}=33$ for any curve $C$. We can now use Proposition \ref{paubry} to
get the bound $q > (2 \cdot 33)^2$.

\subsection{Proof of Theorem \ref{conditionetoile2}: the case $p = 2$}

We first need  to compute the arithmetic genus of $X_C$. It is not  so easy in this case, as wild ramification occurs. We therefore suggest another point of view, which can actually be used in any characteristic.

\begin{lemma}
The arithmetic genus of $X_C$ is $33$.
\end{lemma}
\begin{proof}
In order to emphasize how general the method is, we will denote by $g=3$ the genus of $C$ and by $d=4$ its degree. The map $T$ is
a correspondence with \emph{valence} $\nu=2$, \textit{i.e.} the linear equivalence
class of $T(P)+\nu P$ is independent of $P$. Let denote by $E$ (resp.
$F$) a fiber of $\pi_1$ (resp. $\pi_2$) and $\Delta$ the diagonal of $C \times C$. We get, as in the proof of
\cite[p.285]{griffiths-harris}, that $X_C$ is linearly equivalent to $ a E+b F- \nu \Delta$ for some $a,b
\in \ZZ$  to be determined. 
Then, one computes the arithmetic genus $\pi_{X_C}$ of $X_C$ thanks to the adjunction
formula \cite[Ex.V.1.3.a]{Hart} $$2 \pi_{X_C}-2=X_C.(X_C+\kappa_{C \times
C})$$
where $\kappa_{C \times C}$ is the canonical divisor on $C \times C$.
 Using that (see for instance
\cite[p.288]{griffiths-harris}, \cite[Ex.V.1.6]{Hart})
$$
\begin{cases}
E.E=F.F &=0, \\
E.F=\Delta.E=\Delta.F &= 1,\\
\kappa_{C \times C} & \equiv_{\textrm{num}} (2g-2)E+(2g-2)F,\\
\Delta^2 &= (2-2g),
\end{cases}
$$
 we find
$$\pi_{X_C}=ab-15.$$
Now, we  determine the values of $a$ and $b$. One has  $$X_C.E=b-\nu=\deg \pi_1=d-2=2$$ so $b=4$. Also, 
$X_C.F= a-\nu=\deg \pi_2$. The degree of $\pi_2$ is equal to the degree
of the dual map $C \to C^*$, where $C^*$ is the dual curve, which is $d(d-1)$. In the case $d=4$, we then get $a=12$. Plugging the values of $a$ and $b$, we find  that
$\pi_{X_C}=33$. \qed
\end{proof}

\begin{remark}
In a previous version of this article, we claimed that $a=6$ when $p=2$. This was due to a confusion between the degree of the dual curve and the degree of the dual map, which is in this case inseparable.
\end{remark}

\noindent
We can now come back to the proof of Theorem \ref{conditionetoile2}.\\
 If $X_C$ is absolutely irreducible, then we can argue as in Section \ref{sub:proof} and get the bound $q > (2 \cdot 33)^2$.\\
If $X_C$ is not absolutely irreducible, since $\deg \pi_1=2$, $X_C$ is the union over $\FF_{q^2}$ of two absolutely irreducible curves $X_1,X_2$ such that
$X_1$ and $X_2$ are birationally equivalent to $C$.

Since $C$ is smooth, there are two $\FF_{q^2}$-rational morphisms  $s_i : C \to X_i$ such that $\pi_1|_{X_i} \circ s_i=1$. Hence $\deg s_i=1$ and the $X_i$ are geometrically isomorphic to $C$ and smooth. If the $X_i$ are defined over $\FF_q$ they have a rational point because $q>32$ and $g_{X_i}=3$. This concludes the proof in this case.\\
It is then sufficient to prove that the $X_i$ are defined over $\FF_q$. Otherwise, we would have $X_2=X_1^{\sigma}$, where $\sigma$ is the Frobenius automorphism of $\Gal(\FF_{q^2}/\FF_q)$. Let us consider the morphisms $t_i=\pi_2|_{X_i} \circ s_i : C \to C$ which map a point $P \in C$ to one of the two points in the support of $T(P)$. Since $g_C >1$, these morphisms have to be purely inseparable. Let $P$ be the generic point on $C$ and $Q_1,Q_2=Q_1^{\sigma}$ the two points in the support of $T(P)$. Since $\deg \pi_2=10=\deg \pi_2|_{X_1}+\deg \pi_2|_{X_2}$ and $\pi_2$ is purely inseparable, we have, say, $\deg \pi_2|_{X_1}=2$ and  $\deg \pi_2|_{X_1}=8$. We then get the following commutative diagram

$$\xymatrix{
C \ar[r]^{s_1} \ar[dr]_{s_2} & X_1 \ar[r]^{\pi_2|_{X_1}}   \ar[d]^{F_q} & C \ar[d]^{F_q} \\
& X_2 \ar[r]_{\pi_2|_{X_2}}  & C
}
 \quad \quad 
\xymatrix{
P \ar[r]^{1} \ar[dr]_{1} & (P,Q_1) \ar[r]^{2}   \ar[d]^q & Q_1 \ar[d]^q \\
& (P,Q_2) \ar[r]_{8}  & Q_2
}
$$

where the labelling indicates the degree of the maps and $F_q$ is the $\FF_q$-Frobenius morphism. We hence get a contradiction.

\begin{remark}
Using the first line of the diagram, we see that if $X_C$ is not absolutely irreducible, then there is a purely inseparable degree $2$ morphism from $C$ to $C$. This implies that there is a geometric isomorphism $\gamma :  C^{(2)} \to C$, where $C^{(2)}$ is the conjugate of $C$ by the Frobenius automorphism of $\Gal(\FF_q/\FF_2)$. Hence, the moduli point corresponding to $C$ is defined over $\FF_2$. As the field of moduli is a field of definition over finite fields, $C$ descends over $\FF_2$, \textit{i.e.} there exists  a curve $\tilde{C}/\FF_2$ and a geometric isomorphism $\psi : \tilde{C} \to C$. More specifically, we have the following commutative diagram
$$\xymatrix@C=10mm@R=4mm{
 & C \ar[rr]^{t_1} \ar[dr]_{F_2} & & C \ar[ddr]^{\psi^{-1}} &  \\
 & & C^{(2)} \ar[ur]_{\gamma} \ar[drr]_{(\psi^{-1})^{(2)}} &  & \\
\tilde{C} \ar[uur]^{\psi} & & & &  \tilde{C}
}
$$

Hence the map $\psi^{-1} \circ t_1 \circ \psi : \tilde{C} \to \tilde{C}$ is the $\FF_2$-Frobenius morphism $F_2$. Therefore, the Frobenius morphism maps a point $P \in \tilde{C}$ to a point on the tangent at $P$. Such a (plane) curve $\tilde{C}$ is called Frobenius non-classical. Using \cite[Th.1]{hefez-voloch}, we see that $\tilde{C}$ has no rational point. It is then easy to run a program on plane smooth quartics over $\FF_2$ with no rational points and to check if $X_C$ is irreducible or not. See for instance 
{\small
\begin{center}
\url{http://iml.univ-mrs.fr/~ritzenth/programme/tangent-char2.mag}.
\end{center}
}
 Surprisingly, we find that $X_C$ is reducible if and only if $C/\FF_2$ is isomorphic to 
$$x^4 + x^2 y^2 + x^2 y z + x^2 z^2 + x y^2 z + x y z^2 + y^4 + y^2 z^2 + z^4$$ or
$$x^4 + x^3 z + x y^3 + x y z^2 + y^4 + y z^3 + z^4.$$
These two curves are twists of the Klein quartic. Hence, the curves $C/\FF_q$ such that $X_C$ is reducible are the ones isomorphic to the Klein quartic.
\end{remark}

%% file: Flexes_char_3.tex
\section{The case of flexes} \label{flexcar3}  \label{open}

The heuristic results and computations of \cite{FOR} tend to suggest that a random plane smooth quartic over $\FF_q$ has an asymptotic probability of about $0.63$ to have at least one rational flex when $q$ tends to infinity. We describe how to turn the heuristic strategy into a proof.\\

Let $\PP^{14}$ be the linear system of all plane quartic curves over a field $K$ and $I_0=\{(P,l), \; P \in l\} \subset \PP^2 \times {\PP^2}^*$. Let $I_4 \subset \PP^{14} \times I_0$ be the locus
$$I_4=\{(C,(P,l)), \; C \; \textrm{is smooth and} \; P \; \textrm{is a flex of} \, C\, \textrm{with tangent line} \; l\}.$$   

 Harris proved the following result using monodromy arguments.
\begin{theorem}[{\cite[p.698]{harris-enu}}]
The Galois group of the cover $I_4 \to \PP^{14}$ over $\CC$ is the full symmetric group $S_{24}$.
\end{theorem}
\noindent Let us assume for a moment that this result is still valid over finite fields. Then, using a general Chebotarev density theorem for function fields like in \cite[Th.7]{serre-chebo}, this would mean that the probability of finding a rational flex is within $O(1/\sqrt{q})$
of the probability that a random permutation 
of $24$ letters has a fixed point, which is 
$$p_{24}:=1-\frac{1}{2!} + \frac{1}{3!} - \ldots -\frac{1}{24!} \approx 1-\exp(-1) \approx 0.63.$$
Unfortunately, it is not easy to transpose Harris's proof over any field. And actually, Harris's result is not true in characteristic $3$, as the following proposition implies.

\begin{proposition}
Let $C\,:\, f(x_1,x_2,x_3)=0$ be a smooth plane quartic defined over an algebraically closed field $K$ of characteristic $3$. The flexes of $C$ are  the intersection points of $C$ with a certain curve $H_C : h_C=0$ of degree less than or equal to $2$. The curve $H_C$ can be degenerate as in the case of the Fermat quartic where $h_C = 0$.
\end{proposition}

\begin{proof}
We use the method to compute flexes of a plane curve of degree $d$ described  in the appendix of \cite{FOR} (see also \cite[Th.0.1]{stvo2}). Indeed, when the characteristic of the field divides $2(d-1)$, one cannot use the usual Hessian and one should proceed  as follows.\\
Let  $C : f=0$ be the generic plane quartic over $K$
\begin{eqnarray}
f(x,y,z) &:=& a_{00}y^4+y^3(a_{10}x+a_{01}z)+y^2(a_{20}x^2+a_{11}xz+a_{02}z^2) \nonumber\\
	& & +y(a_{30}x^3+a_{21}x^2z+a_{12}xz^2+a_{03}z^3) \nonumber\\
	& & +(a_{40}x^4+a_{31}x^3z+a_{22}x^2z^2+a_{13}xz^3+a_{04}z^4),	\nonumber
\end{eqnarray}
We define as in \cite{FOR} 
\begin{eqnarray*}
2 \bar{h} &= & 2f_1 f_2 f_{12} - f_1^2 f_{22}-f_2^2 f_{11}, \\
&= &  f_1 ( f_2  f_{12} - f_1 f_{22})+ f_2 ( f_1 f_{12} - f_2 f_{11})
\end{eqnarray*}
where $f_i$ or $f_{ij}$ are the partial derivatives with respect to $i$th variable (or to $i$th and $j$th variables).
It is then easy to check  via a computer algebra system, see

{\small 
\begin{center}
\url{http://iml.univ-mrs.fr/~ritzenth/programme/flex-char3.mw},
\end{center}
}
\noindent that
$$
2\bar{h} -a_{20}f^2-f \cdot  (ax^3+by^3+cz^3)z = \tilde{h}_{C} \cdot z^2, $$
with
\begin{eqnarray}
a &:= & a_{40}a_{11}+a_{21}a_{30}-2a_{20}a_{31},  \nonumber\\
b	&:= & a_{10}a_{11}+a_{00}a_{21}-2a_{20}a_{01},\nonumber \\
c	&:= & 2a_{12}^2+a_{13}a_{11}+a_{20}a_{04}+a_{03}a_{21}+a_{02}a_{22},\nonumber
\end{eqnarray}
where $\tilde{h}_{C}$ is a homogeneous polynomial in $K\left[x,y,z\right]$ of degree $6$, for which nonzero coefficients appear  only for the monomials $ x^6,y^6,z^6,x^3y^3,x^3z^3$ and $y^3z^3$.
  Since the map $u\mapsto u^3$ is an isomorphism of $K$, there is a polynomial $h_{C}\in K\left[x,y,z\right]$ satisfying $\tilde{h}_{C}=h_{C}^3$. If we suppose that there is no flex at infinity, the flexes are the intersection points of  $\bar{h}=0$ and $f=0$, so they are also the intersection points of $h_C=0$ and $f=0$ and $h_C=0$ is the equation of a (possibly degenerate) conic $H_C$. \qed
\end{proof}
\begin{remark}
As $\tilde{h}_C=h_C^3$, we see that the weight of a flex which is not a hyperflex is $3$ in characteristic $3$.
\end{remark}

\begin{corollary}
The Galois group of the cover $I_4 \to \PP^{14}$ over $\bar{\FF}_3$ is the full symmetric group $S_{8}$.
\end{corollary}
\begin{proof}
Note that the Galois group $G$ of the cover is the Galois group of the $x$-coordinate of the $8$ flexes of the general quartic. Hence, $G$ is included in $S_8$. To show that $G$ is exactly $S_8$, we are going to specialize the general quartic to smooth quartics over finite fields with $8$ distinct flexes having different arithmetic patterns. More precisely, 
to generate $S_8$ we need to produce (see \cite[Lem.4.27]{milne-galois}):
\begin{itemize}
\item a smooth quartic with $8$  Galois conjugate flexes over $\FF_3$:
$$2 x^4 + 2 x^3 y + x^3 z + 2 x^2 z^2 + x y^3 + 2 x y^2 z + y^3 z + y z^3=0 ;$$
\item a smooth quartic with one rational flex and $7$ Galois conjugate flexes over $\FF_3$:
$$x^3 y + x^2 z^2 + 2 x y^3 + x y^2 z + 2 x y z^2 + 2 x z^3 + y^4 + 2 y z^3=0;$$
\item a smooth quartic with two quadratic conjugate flexes and $6$ rational flexes over $\FF_9$: 
\begin{eqnarray*}
a^6 x^4 &+& a x^3 y  +  a^7 x^3 z + a^6 x^2 y^2 + a^2 x^2 z^2 + a^7 x y^3 + a^7 x y^2 z \\ &+& x y z^2 
 +  a^5 x z^3 + a^5 y^4 + a^3 y^3 z + a^5 y^2 z^2 + 2 y z^3 + a^7 z^4=0,
\end{eqnarray*}
where $a^2-a-1=0$. \qed
\end{itemize} 

\end{proof}
\begin{corollary}
Let $C$ be a smooth plane quartic over $\FF_{3^n}$. The probability that $C$ has a rational flex tends to 
$$p_8:=1-\frac{1}{2!} + \frac{1}{3!} - \ldots -\frac{1}{8!}  \approx 0.63$$
when $n$ tends to infinity.
\end{corollary}
\begin{remark}
The fact that the Galois group $G$ in characteristic $3$ is $S_8$ and not $S_{24}$ was unnoticed in our computations in \cite{FOR} because $|p_{24}-p_8| \leq 10^{-5}$.
\end{remark} 
General reduction arguments show that the Galois group $G$ remains $S_{24}$  almost all $p$. 
We  conjecture that $p=3$ is the only exceptional case.
\begin{conjecture}
The Galois group of the cover $I_4 \to \PP^{14}$ over $\bar{\FF}_p$ is the full symmetric group $S_{24}$ if $p \ne 3$ and $S_8$ otherwise.
\end{conjecture}